\documentclass[11pt,a4paper]{article}
\usepackage[utf8]{inputenc} 

\usepackage{epsfig}
\usepackage{epstopdf}
\usepackage{graphicx}
\graphicspath{{Figures/}}
\usepackage{subfig}
\usepackage{color}
\usepackage{pseudocode}
\usepackage{url}

\usepackage{amssymb}
\usepackage{amsmath}
\usepackage{amsthm}

\newtheoremstyle{meiner} 
    {4pt}{4pt}           
    {\itshape}           
    {}                   
    {\sffamily\bfseries} 
    {}                   
    { }                  
    {}                   
\theoremstyle{meiner}

\newtheorem{thm}{Theorem}[section]
\newtheorem{lem}[thm]{Lemma}
\newtheorem{prop}[thm]{Proposition}

\theoremstyle{definition}

\theoremstyle{remark}

\def\th{{\ensuremath{{}^{\text{\tiny th}}}}\ }
\def\Proof{\par{\NI\bf Proof. }}
\def\qed{\hfill\fbox{\hbox{}}\bigskip}

\def\Case#1.{\rih{Case #1.}}
\def\Phase#1.{\NI{\bfPhase #1.}}
\def\Fact#1.{\par\smallskip{\NI\bf Fact #1.}}
\def\Claim#1.{{\bf Claim #1.}}
\def\ClaimX{\par\smallskip{\NI\bf Claim.}}

\def\forget#1{}

\def\NI{\noindent}
\def\ni{\noindent}

\def\ITEMMACRO #1 ??? #2 ???{\par\medskip\noindent%
\hangindent=#2em\setbox0\hbox{#1 \kern5pt}%
\ifdim\wd0<\hangindent\setbox0\hbox to\hangindent{\hss#1\quad}\fi%
\box0\ignorespaces}

\def\Item#1{\ITEMMACRO #1 ??? 2.5 ???}

\def\Bitem{\Item{\hss$\bullet$}}

\def\abstract{\noindent\hfil\vbox\bgroup\hsize=.9\hsize
\small\noindent{\bf Abstract.}
}
\def\endabstract{\egroup\hfil}

\def\width{{\sf width}}
\def\sp{{\sf spread}}
\def\rk{{\sf rk}}
\def\PsFig#1#2#3{}
\def\Pr{{\sf Pr}}

\begin{document}

\title%
 {Linear Extensions of N-free Orders}

\author{
\parbox{6.5cm}{\center
{\sc Stefan Felsner
\footnote{Partially supported by DFG grant FE-340/7-2 and 
              ESF EuroGIGA project ComPoSe.}}\\[3pt]
\normalsize
\url{felsner@math.tu-berlin.de}\\
\small
        {Institut f\"ur Mathematik\\
         Technische Universit\"at Berlin\\
         Strasse des 17. Juni 136\\
         D-10623 Berlin, Germany}}
\and
\parbox{6.5cm}{\center
{\sc Thibault Manneville}\\[3pt]
\normalsize
\url{tmannevi@clipper.ens.fr}\\
\small
      {{\'E}cole Normale Sup{\'e}rieure de Paris\\ 
       45, rue d'Ulm\\
       75005 Paris, France}}
}

\date{}

\maketitle

\abstract We consider the number of linear extensions of
an N-free order $P$. We give upper and lower bounds on this number in
terms of parameters of the corresponding arc diagram. We propose a
dynamic programming algorithm to calculate the number. The algorithm
is polynomial if a new parameter called activity is bounded by a
constant. The activity can be bounded in terms of parameters of the
arc diagram. 
\endabstract
\bigskip

\textbf{Mathematics Subject Classifications (2000)} 06A05, 06A07, 68R99 
\bigskip

\textbf{Keywords :} arc diagram, counting linear extensions, bounds for the number of linear extensions, dynamic programming 

\section{Introduction}

The number of linear extensions is one of the most
fundamental combinatorial parameters of an order
(poset).  Explicit formulas and efficient algorithms have been found
for several classes of orders, see e.g.~\cite{atkinson-bounded-width}
\cite{linear-extensions-of-trees},
\cite{elementary-proof-of-the-hook-formula} or
\cite{Habib-Mohring-qsp-posets}.  Brightwell and
Winkler~\cite{P-completude} have shown that for general orders the
problem is \#P-complete. Hence, unless P=NP, there is no polynomial
algorithm for this problem.

Problems that are hard for general orders may become tractable when
restricted to more structured classes of orders,
see~\cite{computationally-tractables-cases} for a survey on the topic. One of the most
prominent classes of orders is the class of N-free ordered sets. This
class was introduced by Grillet~\cite{grillet-max-chains-antichains}.
Leclerc and Montjardet~\cite{ordres-cac-leclerc} characterized N-free
orders as the chain-antichains-complete orders. A series of papers
investigated algorithmic problems on N-free orders : Habib and
Jegou~\cite{N-free-posets-as-generalization-of-sp} and independently
Sys\l{}o~\cite{minimizing-jump-number-syslo} showed that the jump
number of an N-free order can be computed in linear time. For general
orders the jump-number problem is NP-hard. Habib and
M\"ohring~\cite{Habib-Mohring-qsp-posets} proved that the isomorphism
problem remains isomorphism-complete when restricted to N-free
orders. More recent contributions to N-free orders are the enumeration
of N-free orders of a given
height~\cite{height-counting-interval-and-N-free-posets} and the
computation of their page number~\cite{page-number-of-N-free-posets}.
Zaguia~\cite{1-2-2-3-conjecture-N-free-posets} verified the $1/3-2/3$
conjecture for N-free orders.
We are not aware of previous results concerning the number of linear
extensions of N-free orders.

The next section provides definitions and background information about
N-free orders.  In Section~\ref{sec:bounds} we present upper and lower
bounds for the number of linear extensions of an N-free order.

Section~\ref{sec:algo} is devoted to the dynamic programming algorithm
for calculating the number of linear extensions of N-free orders. 
Dynamic programming was used to show that the number of linear extensions
of orders of bounded width can be computed in polynomial time, see~\cite{atkinson-bounded-width}.
In our approach the size of an antichain is
replaced by the size of an {\it active set}. In particular our algorithm
has a polynomial running time for some N-free orders of unbounded
width.

We conclude with a discussion of some open problems.

\subsection{Definitions and background on N-free orders}

We are concerned here with combinatorial problems for partially ordered sets.
We assume some familiarity with concepts and results in this area, including
linear extensions, comparability graphs and diagrams.  For readers who are new
to the subject, we suggest consulting one of the books on the topic
(Trotter~\cite{trotter-dimension-theory}, Schr\"oder~\cite{schroeder-ordered-sets}, Neggers and
Kim~\cite{negger-kim-basic-posets}).

An order is \emph{N-free} if its diagram does not contain an N (see
Figure~\ref{N poset}).  There are many characterizations for N-free orders, an
overview is given by M\"ohring~\cite{computationally-tractables-cases}.  The relevant
characterization in our context is the following :

\begin{figure}[hbtp]
 \begin{center}
  \ \psfig{file=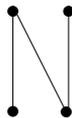,height=1.5cm}
 \end{center}
 \caption{The order N.} \label{N poset}
\end{figure}

\begin{thm}\label{caracterisation des N-free}
An order is N-free if and only if its diagram is the line digraph 
of a directed acyclic graph (dag).
\end{thm}

\begin{figure}[hbtp]
 \begin{center}
  \ \psfig{file=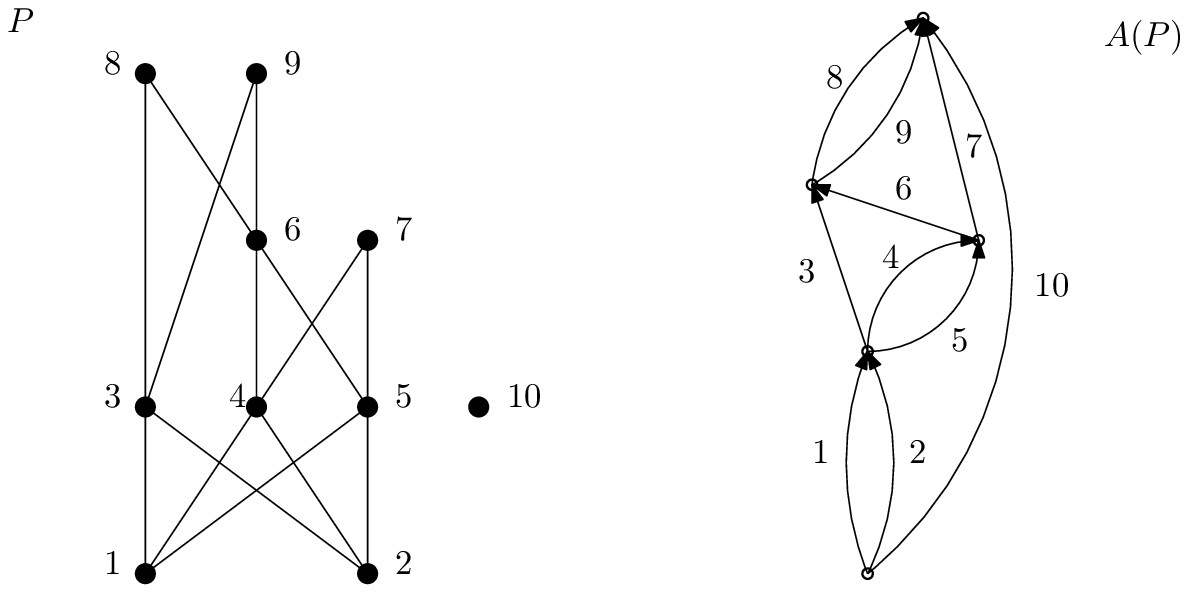}
 \end{center}
 \caption{An N-free order and its arc diagram.} \label{fig:arc-diagram}
\end{figure}

Theorem \ref{caracterisation des N-free} enables us to represent
N-free orders, using an \emph{arc diagram}. The arc diagram $A(P)$ of
an N-free order $P$ is the digraph of which the order is the line
digraph. Figure \ref{fig:arc-diagram} shows an example. Such a digraph is
not unique, but if we require that $A(P)$ has a unique source $\hat{0}$ and a
unique sink $\hat{1}$, then it becomes unique. 
Sys\l{}o~\cite{minimizing-jump-number-syslo}
used $A(P)$ to deal with jump number of N-free
orders. Since $A(P)$ can easily be computed we will assume that $P$ and $A(P)$ 
are both available. Depending on the context we then think of elements of $P$ 
as vertices of $P$ or as edges of $A(P)$, whatever is more convenient.

Being acyclic the arc diagram $A(P)$ or more precisely its transitive closure
is a partial order relation on the vertex set of $A(P)$.  In the following, we
will switch hence and forth between the arc diagram as a directed graph
and the arc diagram as an order. The interpretation will be clear from the
context. It is important to notice that in general the digraph $A(P)$ is not
the diagram of the order $A(P)$, simply because $A(P)$ may have transitive
edges.

\section{Bounds for the linear extensions of an N-free order}\label{sec:bounds}

Using the arc diagram representation of an N-free order, we prove the following

\begin{prop}
 Let $P$ be an N-free order and $e(P)$ be its number of linear extensions, then 
 $$
e(A(P))\prod_{v\in V_{A(P)}}o_v!
\quad\leq\quad
e(P)
\quad\leq\quad
|P|!\prod_{v\in V_{A(P)}}\binom{i_v+o_v}{o_v}^{-1}
$$
 where $V_{A(P)}$ is the set of the vertices of $A(P)$, $o_v$ is the
 out-degree of $v\in V_{A(P)}$, and $i_v$ is the in-degree of $v\in V_{A(P)}$.
\end{prop}

\Proof For the lower bound we construct an appropriate set of linear
extensions of $P$. Let $L=(v_1,v_2,..,v_k)$ be a linear extension (topological
ordering) of the arc diagram $A(P)$.  For each vertex~$v$ of $A(P)$, choose a
permutation $\sigma_v=(x_1^v,x_2^v,..,x_{o_v}^v)$ of the out-edges of $v$.
Concatenating these permutations in the order given by $L$ yields
$L'=(x_1^{v_1},x_2^{v_1},..,x_{o_v}^{v_1},\ldots,x_1^{v_i},..,x_{o_{v_i}}^{v_i},
\ldots,x_{o_{v_k}}^{v_k})$, a permutation of the elements of $P$. We claim that
$L'$ is a linear extension of $P$. If $x < y$ in $P$, then there is a
directed path in $A(P)$ whose first and last edge are $x$ and $y$
respectively. Let $v_x$ and~$v_y$ be the source vertices of the edges $x$ and
$y$. The directed path implies the relation $v_x < v_y$ in $A(P)$, whence
$v_x$ precedes $v_y$ in $L$. This in turn implies that $x$ precedes $y$ in
$L'$. We have thus shown that $L'$ respects the order relations of $P$, i.e.,
it is a linear extension.

Hence, there is a linear extension $L'$ for every tuple
$(L,\sigma_1,\ldots,\sigma_k)$. Obviously, if $L'$ and $L''$
are constructed using different tuples, then they are different.
The lower bound is just the number of tuples.

The upper bound is based on a probabilistic argument.  Let $\sigma$ be a
randomly chosen permutation of the elements of $P$. The upper bound will be
obtained by considering the probability that~$\sigma$ is a linear
extension of $P$.

A vertex $v\in A(P)$ is \emph{good} with respect to $\sigma$, in symbols
$v\!\blacktriangleright\!\sigma$, if in $\sigma$ all in-edges of $v$ 
precede all out-edges of $v$. The probability that $v$ is
good with respect to a random permutation~$\sigma$ is $\binom{i_v+o_v}{o_v}^{-1}$:
just consider the induced permutation $\sigma_v$ on the edges incident to $v$,
out of the $(i_v+o_v)!$ possible permutations~$\sigma_v$ exactly $i_v!o_v!$ 
make $v$ good.

Now we try to see what happens if we have to deal with more than one
vertex. Let $\Pr(A)$ denote the probability of $A$ and
$\Pr(A|B) = \Pr(A \cap B)\Pr(B)^{-1}$ be the conditional probability of $A$
given $B$.
Since a permutation is a linear extension of $P$ if and only if all
the vertices of $A(P)$ are good for this permutation, we get:
\begin{equation} \label{proba upper bound N-free}
  \Pr(\sigma\in\mathcal{L}(P))=
  \prod_{j=1}^{k} \Pr(v_{j}\!\blacktriangleright\!\sigma|
           v_1\!\blacktriangleright\!\sigma,\ldots v_{j-1}\!\blacktriangleright\!\sigma)
\end{equation}
Where $(v_1,\ldots,v_k)$ is an arbitrary enumeration of the vertices of $A(P)$
and $\mathcal{L}(P)$ is the set of the linear extensions of $P$. To shorten
the proof a little, we will assume that $(v_1,\ldots,v_k)$ is a linear
extension of $A(P)$.

We claim that for all $j\in[k]$
\begin{equation}\label{correlation}
\Pr(v_j\!\blacktriangleright\!\sigma|v_1\!\blacktriangleright\!\sigma,\ldots,
v_{j-1}\!\blacktriangleright\!\sigma)
\leq
\Pr(v_{j}\!\blacktriangleright\!\sigma)=\binom{i_{v_j}+o_{v_j}}{o_{v_j}}^{-1}.
\end{equation}
Plugging these inequalities into (\ref{proba upper bound N-free}) yields the 
upper bound.

Permutations $\sigma$ and $\pi$ of $P$ are called {\it equivalent for $j$} if: 
\Item{1.} all $x\in P$ that are not incident to $v_j$ are at the same position in
$\sigma$ and $\pi$, 
\Item{2.} both permutations induce the same orderings (permutations) on the in-edges and the
out-edges of $v_j$, i.e., 
$\sigma_j^{{\sf in}}=\pi_j^{{\sf in}}$ and 
$\sigma_j^{{\sf out}}=\pi_j^{{\sf out}}$.

\medskip

\ni
\ClaimX\ If an equivalence class $\cal C$ contains
a permutation that makes $v_j$ good, then\\
\centerline{$\Pr(v_j\!\blacktriangleright\!\sigma|\sigma\in {\cal C}) = \binom{i_{v_j}+o_{v_j}}{o_{v_j}}^{-1}$,
otherwise $\quad\Pr(v_j\!\blacktriangleright\!\sigma|\sigma\in {\cal C}) =0$.}
\medskip

To see that the claim implies the correlation inequality
(\ref{correlation}) we propose an experiment in two phases. First we pick
a class with a probability proportional to its size, then we pick a
random element from the class. When the probability space for the
experiment consists of all permutations, the probability for picking a $\sigma$
with $v_j\!\blacktriangleright\!\sigma$ is
$\binom{i_{v_j}+o_{v_j}}{o_{v_j}}^{-1}$.  For a proof of (\ref{correlation})
we are interested in a probability space
for the experiment that only consists of permutations satisfying
$v_1\!\blacktriangleright\!\sigma,\ldots,
v_{j-1}\!\blacktriangleright\!\sigma$. In this space we only get a $\sigma$ with
$v_j\!\blacktriangleright\!\sigma$ if we are lucky in the first phase
of the experiment. From the claim we know that, conditioned on success
in the first phase, the probability of success in the second phase
is again $\binom{i_{v_j}+o_{v_j}}{o_{v_j}}^{-1}$. This yields~(\ref{correlation}).
\medskip

\ni {\bf Proof of the claim.} Let $\cal C$ be a class. With a permutation
$\sigma \in\cal C$ consider the induced permutation $\sigma_j$ on the edges
incident to $v_j$. From the two properties in the definition of equivalence it
follows that $\sigma$ is uniquely determined if we know which positions of
$\sigma_j$ are used by in-edges and which by out-edges. This can be encoded by
a 0-1 vector $s_\sigma=(s_1,s_2,\ldots,s_{i_{v_j}+o_{v_j}})$ with $s_i = 1$ if
the $i$\th entry of $\sigma_j$ is an in-edge and $s_i = 0$ if it is an
out-edge. 

Now consider an adjacent inversion $(s_i,s_{i+1})=(1,0)$ in $s_\sigma$ and let
$\sigma'$ be the permutation obtained from $\sigma$ by swapping the out-edge
$x$ and the in-edge $y$ of $v_j$ that are at positions $i$ and $i+1$ in
$\sigma_j$.  To show that $\sigma'$ belongs to the same class it suffices to
show that $v_\ell\!\blacktriangleright\!\sigma'$ for all $1\leq\ell<j$.  Let
$x$ be the edge $v_j \to v_{j_x}$ and $y$ be the edge $v_{j_y}\to v_j$.  Since
$(v_1,\ldots,v_k)$ is a topological order of $A(P)$ we have $j_y<j<j_x$. We
don't care whether $v_{j_x}$ is good with respect to $\sigma'$.  Vertex
$v_{j_y}$ is good for $\sigma$ and in $\sigma'$ an out-edge of $v_{j_y}$ has
been moved further to the right compared to $\sigma$, hence, $v_{j_y}$ is good
for $\sigma'$. For all the other $\ell<j$ the property
$v_\ell\!\blacktriangleright\!\sigma'$ is directly inherited from
$\sigma$. Note that the encoding vector $s_{\sigma'}$ of $\sigma'$ has exactly
one inversion pair less than~$s_{\sigma}$.

Let $s$ and $s'$ be two 0-1 vectors of length $t$ with $r$ ones such that $s'$
is majorized by $s$, i.e., $\sum_{i=1}^j s'_i \leq \sum_{i=1}^j s_i$ for all
$j=0,..,t$. Sequence $s'$ can be reached from $s$ by a sequence of steps that
remove a single inversion each. Now, suppose class $\cal C$ contains a
permutation $\sigma$ that makes $v_j$ good. Being good is equivalent to
$s_\sigma = (1,..,1,0,..,0)$ which is the maximum in the majorization
order. Hence, starting from $\sigma$ we can find a sequence of inversion
reducing swaps to reach every permutation $\sigma'$ such that $s_{\sigma'}$ is
majorized by $s_\sigma$. Since $\sigma\in \cal C$ and the steps do not leave
$\cal C$ we find that $|{\cal C}|=\binom{i_{v_j}+o_{v_j}}{o_{v_j}}$.  This
completes the proof of the claim and of the proposition.~\qed

\begin{figure}[hbtp]
 \begin{center}
  \ \psfig{file=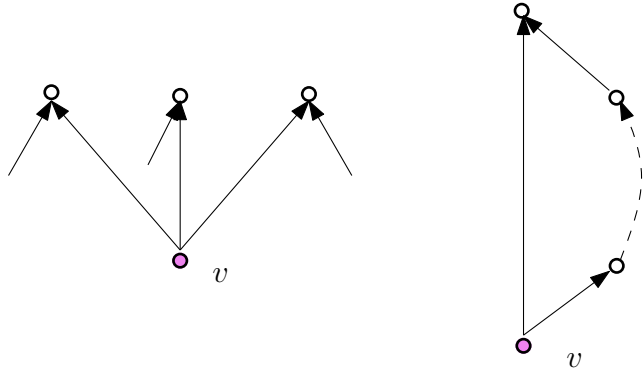}
 \end{center}
\caption{Subdigraphs of $A(P)$ that allow the addition of the lower bound and its dual.} \label{fig:subdigs}
\end{figure}

There are few examples of N-free orders where one of the bounds is
sharp.  In special cases, however, there may be many ways for
improving the bounds.  We exemplify this with the lower bound.  First
note that the dual $P^*$ of an N-free order $P$ is N-free again and
$A(P^*)$ is the dual of $A(P)$.  Indeed $A(P^*)=A(P)^*$ holds in both
interpretations, as a digraph and as an order. Hence the dual lower
bound can be written as $e(A(P))\prod_{v\in V_{A(P)}}i_v!$ and also
provides a lower bound on $e(P)$. Moreover, in many cases the set of
linear extensions of $P$ that are obtained from the lower bound proof
and the dual are disjoint. This is always the case if there is a
vertex $v$ in $A(P)$ such that the out-edges of $v$ are always
separated if we keep the in-edges of all vertices together. In
particular, if $A(P)$ has a subdigraph of one of the types shown in
Figure~\ref{fig:subdigs} then $e(A(P))(\prod_{v\in
  V_{A(P)}}o_v!+\prod_{v\in V_{A(P)}}i_v!)$ is a lower bound on $e(P)$.

\section{A dynamic programming algorithm for enumeration}\label{sec:algo}

In this section we describe a dynamic programming algorithm for the
enumeration of linear extensions of an N-free order.

Let $P$ be an N-free order with arc diagram $A(P)$.  Given a linear
extension $(x_1,\ldots,x_n)$ of $P$ the algorithm computes the number
of linear extensions of the induced suborders $P_i = P[x_1,..,x_i]$ for 
$i$ increasing from $1$ to $n$. Denote the edge of $A(P)$ that represents $x_j$ by
$(v_j^-,v_j^+)$.  The \emph{active set} $A_i$ corresponding to $P_i$
consists of all vertices of $A(P)$ that have incident edges from $P_i$
and from the complement of $P_i$, i.e., $A_i = \{ v \;:\; \exists\; j'
\leq i < j'' \text{\;with\;} v_{j'}^+ = v = v_{j''}^- \}$.  Let $L$ be
a linear extension of $P_{i-1}$, extending $L$ to a linear extension of
$P_i$ means that $x_i$ is placed at a position behind all of its
immediate predecessors in $L$.  The immediate predecessors of $x_i$
are exactly the elements that point into $v_i^-$, whence $v_i^-$
belongs to the active set $A_{i-1}$. 

We are now ready to describe the algorithm.
For $i=1,\ldots,n$ and a function
$f:A_i\to[i]$ we define
\begin{eqnarray*}
z_i(f) & =\;\; \#\big(\;L \;:\!\!\!&\text{linear extension of\;}P_i\text{\;such that in\;} L 
         \text{\;the last element \;} x_j  \text{\;with\;} v_{j}^+=v\\
       & &  \text{\;is at position\;} f(v), \text{\;for all\;} v \in A_i).
\end{eqnarray*}
Let $z_{i-1}(f)$ be known for all $f:A_{i-1}\to[i-1]$, then we can compute 
$z_i(f)$ for all $f:A_{i}\to[i]$.  The first thing that has to be done is to
find $A_i$. Only two events have to be checked for the update from $A_{i-1}$ 
to $A_i$ :
\Bitem If all element incident to $v_i^-$ belong to $P_i$, then $v_i^-$ is removed from the active set.
\Bitem If $x_i$ is the only element of $P_i$ incident to $v_i^+$, then $v_i^+$ is added to the active set. 
\smallskip

\ni 
For fixed $f$ let consider a linear extension $L$ counted by
$z_{i-1}(f)$. Inserting $x_i$ in $L$ after the $t$\th element yields a
linear extension of $P_i$ if and only if $f(v_i^-) \leq t < i$. If $t$
respects these bounds we obtain a linear extension of $P_i$ counted by
$z_i(f')$ where for $v\in A_i$ the value $f'(v)$ is defined as follows:
$$
f'(v) = \begin{cases} 
           \max(f(v),t)+1  &\mbox{if } v=v_i^+\\
           f(v)    &\mbox{if } f(v) \leq t\\
           f(v)+1  &\mbox{if } f(v) >  t.
        \end{cases}
$$
We write $f'=\Gamma_i(f,t)$. The initialization for the dynamic
program is $z_i(f)=0$ for all $i$ and $f$ except $z_1(f_0) = 1$, where
$f_0(v_1^+) = 1$ (note that $A_1=\{v_1^+\}$). The outer loop of the
algorithm is increasing~$i$ from~$2$ to $n$. For each $i$ we go
through all $f$ and all $t$ with $f(v_i^-)\leq t \leq i$ and update
$z_i(\Gamma_i(f,t)) \gets z_i(\Gamma_i(f,t)) + z_{i-1}(f)$.  

Let
$\alpha$ be an upper bound on the size of the active sets, i.e, if
$|A_i| \leq \alpha$ for all $i$. For each of the~$n$ levels $i$ of the 
computation we have to consider at most $n^\alpha$ functions $f$.
For each function~$f$ we have to consider $\leq n$ derived functions $f'$,
if in the loop for $f$ we decrease $t$ in steps of $-1$, then updates 
can be done in $O(1)$. Altogether this yields a running time of
$O(n^{\alpha+2})$.

We can improve the running time by stripping off a factor of $n$.
The key observation is that in each relevant $f$ at level $i$, i.e., each $f$
with $z_i(f) \neq 0$, we have some $v\in A_i$ with $f(v) = i$.
Hence, if we mark this element $v$ we can describe the rest of 
$f$ as a function in $[n]^{\alpha-1}$.
This shows that we only have to consider $\alpha n^{\alpha-1}$ functions on
each level. 

In the following subsection we discuss conditions that allow
to bound the size of the active set.

\subsection{Bounding the size of active sets}\label{ssec:condi}

Since each active set $A_i$ is a subset of the vertex set of the arc diagram
we easily get:
\begin{prop}\label{prop:v-size}
The number of linear extensions of N-free orders with $n$ elements and 
an arc diagram with $\leq k$ vertices can be computed in $O(n^{k+1})$.
\end{prop}

The initial linear extension $L_0=(x_1,\ldots,x_n)$ of $P$ determines the
sequence $A_1,\ldots,A_n$ of active sets. Different choices for $L_0$ may
lead to active sets of very different size. Let us define 
the \emph{activity $\alpha(L)$ of a linear extension} as the maximum size
of an active set in its sequence. The \emph{activity an N-free order} is
$$
\alpha(P) = \min\big( \alpha(L) \;:\; L \text{\; a linear extension of \;} P \big)
$$
In some sense the strongest result that can be stated on the basis of our
algorithm is the following
\begin{prop}
The number of linear extensions of an N-free order $P$
of bounded activity can be computed efficiently if a linear
extension $L$ with $\alpha(P)=\alpha(L)$ is given.
\end{prop}

However, we do not know how to compute the activity of an N-free
order. We propose the activity as a new parameter of N-free
ordered sets for further studies.

\begin{figure}[hbtp]
 \begin{center}
  \ \psfig{file=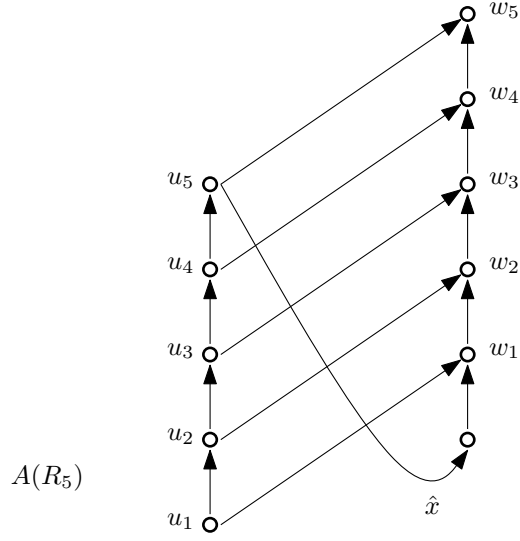}
 \end{center}
\caption{The arc diagram of an N-free order $R_5$ with $\alpha(R_5) = 6$.} \label{fig:counterex}
\end{figure}

It would be nice to improve on Proposition~\ref{prop:v-size} by only bounding
the width of the arc diagram. This, however, is not enough. In
Figure~\ref{fig:counterex} we show the arc diagram $A(R_\ell)$ of an N-free
order $R_\ell$ with $3\ell$ elements and $\width(A(R_\ell))=1$, the activity
of $R_\ell$ is at least $\ell$.  Indeed if $\hat{x}$ is at position~$i$ in
$L$, then $A_i$ contains at least one of $\{u_i,w_i\}$ for $i=1,..,\ell$.

For two vertices $u$ and $v$ of a directed acyclic graph $D$ we define
$\sp(u,v)$ by looking at pairs of interiorly disjoint $u$ to $v$
paths. If there is no such pair, then $\sp(u,v)=0$, otherwise we let
$\sp(u,v)$ be the maximal difference in length of two such paths. The
\emph{spread} of $D$ is $\sp(D) = \max( \sp(u,v) \;:\; u,v\in D )$.
Since shortest and longest paths in directed acyclic graphs are efficiently
computable the spread is also tractable.

Note that $\sp(A(R_\ell)) = \ell+1$. In fact, the activity of an
$N$-free order $P$ is small if width and spread of $A(P)$ are small.

\begin{lem}
For N-free orders $P$ we have
$\alpha(P) \leq \width(A(P))\cdot(\sp(A(P))+2)$.
\end{lem}

\Proof The rank ${\sf rk}(v)$ of a vertex of $A(P)$ is the length of a longest
path ending in $v$. Let $T=(v_1,\ldots,v_s)$ be a topological order of $A(P)$
such that $i < j$ implies ${\rk}(v_i) \leq {\rk}(v_j)$.  Define a linear
extension $L_T$ such that for $i < j$ the out-edges of $v_i$ precede the
out-edges of $v_j$.  

\ClaimX\ $\alpha(L_T) \leq \width(A(P))\cdot(\sp(A(P))+2)$.
\smallskip

Partition the vertices of $A(P)$ into ranks $S_j = \{ v \in V_{A(P)} \;:\;
\rk(v) = j \}$.  By definition each $S_j$ is an antichain of $A(P)$, hence
$|S_j|\leq \width(A(P)$.  If $x_i$, the $i$\th element of $L_T$, has its
source $v_{x_i}^-$ in $S_j$, then $A_i \subset S_j \cup S_{j+1} \cup\ldots\cup
S_{j+\sp(A(P))+1}$. To see this consider $v\in A_i$ with $\rk(v) =k$. Let $p$
be a path of length $k$ from $\hat{0}$ to $v$ in $A(P)$.  From $v\in A_i$ and
the construction of $L_T$ we conclude that there is $y\in P$ with $\rk(v_y^-)
\leq j$ and $v_y^+=v$. Let $q$ be a path from $\hat{0}$ through $v_y^-$ to
$v$. The length of $q$ is at most $j+1$. Note that after clipping $p$ and $q$
from $\hat{0}$ to the last common vertex $u$ we obtain a pair of disjoint $u$
to $v$ paths whose difference in length is at least $k - (j +1)$, hence
$\sp(A(P)) \geq k - j -1$.  \qed

\begin{thm}
The number of linear extensions of N-free orders with $n$ elements and 
an arc diagram with width $\leq k$ and spread $\leq s$ can be computed in 
$O(n^{k(s+2)+1})$.
\end{thm}

\section*{Conclusion}

We have shown that the number of linear extensions of N-free orders of
bounded activity can be computed in polynomial time.  It would be
interesting to understand the class of N-free orders of bounded
activity better. Is membership in this class testable in polynomial
time? Are there simple additional sufficient conditions or necessary
conditions for membership in this class?

We do not expect that there is a polynomial algorithm that solves the
problem for general N-free orders.  Indeed we conjecture that counting
linear extensions of N-free orders is $\#$P-complete.

Brightwell and Winkler~\cite{P-completude} have shown
$\#$P-completeness for general orders of height $3$.  They expect that
the problem remains hard for orders of height~$2$. Since N-free orders
of height $2$ are disjoint unions of weak orders of height $2$ the
counting problem for this class is easy. We think that the counting
problem should already be hard for  N-free orders of height $3$.


\bibliography{bibliography_posets}
\bibliographystyle{my-siam}

\end{document}